\documentclass[12pt]{amsart}   
\usepackage{amsmath,amssymb}   
\usepackage{a4wide} 

\newcommand{\R}{{\mathbb R}}      
     
\newcommand{\N}{{\mathbb N}}      
\newcommand{\T}{{\mathbb T}}      
\newcommand{\Z}{{\mathbb Z}}   

\newcommand{\ve}{\varepsilon}   
    
\newcommand{\Lo}{L^2 (\varOmega,m)}
\newcommand{\pos}{p_\omega (D)}
\newcommand{\Del}{\mathcal{D}_{r,R}}

\DeclareMathOperator{\card}{card}

\theoremstyle{plain} 
\newtheorem{theorem}{Theorem}         
\newtheorem{lemma}{Lemma}         
         
\newtheorem{coro}{Corollary}    

\theoremstyle{definition}   
\newtheorem{definition}{Definition}         
\newtheorem{remark}{Remark}   
\newtheorem{ex}{Example}

\begin{document}      
\title[Pure point diffraction and zero entropy]
{Pure point diffraction implies zero entropy \\[1mm]
 for Delone sets with uniform cluster frequencies}

\author{Michael Baake}
\address{Fakult\"at f\"{u}r Mathematik, Universit\"{a}t Bielefeld, \newline
\hspace*{1.35em}Postfach 100131, 33501 Bielefeld, Germany}
\email{$\{$mbaake,richard$\}$@math.uni-bielefeld.de}

\author{Daniel Lenz}
\address{Fakult\"at f\"ur Mathematik, TU Chemnitz,
09107 Chemnitz, Germany}
\email{dlenz@mathematik.tu-chemnitz.de}

\author{Christoph Richard}

\begin{abstract}   
   Delone sets of finite local complexity in Euclidean space are
   investigated.  We show that such a set has patch counting and
   topological entropy $0$ if it has uniform cluster frequencies 
   and is pure point diffractive.  We also note that the patch
   counting entropy vanishes whenever the repetitivity function 
   satisfies a certain growth restriction.
\end{abstract}      

\maketitle
 

\bigskip
\bigskip
      
\section{Introduction}  \label{Introduction}   

Aperiodic order is a regime of (dis)order at the very border between
order and disorder. It has attracted a lot of attention in recent
years, see \cite{AG95,M97,P95,BM00} and references therein, though
many of its facets are still rather enigmatic.  On the one hand, this
attention originated in intriguing consequences of aperiodic order. On
the other hand, it is due to the existence of real world solids that
exhibit this form of order. These substances are called quasicrystals,
and were discovered (in 1982) by their diffraction properties
\cite{Shecht}.  More precisely, they exhibit pure point diffraction
(indicating long range order), while typically displaying
non-crystallographic symmetries at the same time.

So far, there is no axiomatic framework for aperiodic order, though
various fingerprints of order and their relationships have been
considered.  Besides pure point diffraction, these properties include
finite local complexity (FLC), the Meyer property, repetitivity, and
uniform cluster frequencies (UCF).  In fact, assuming finite local
complexity, it is possible to study the patch counting complexity
(compare \cite{Lag,LP1,LP2}).  A first (though still rather coarse)
measure of order is then given by the patch counting entropy, which
turns out to be finite, compare \cite{Lag,LP1}. Of course, in an
ordered situation, it is expected to vanish.  Thus, both pure point
diffraction and vanishing patch counting entropy can be seen as
indications of an underlying order, and it is natural to ask for their
relationship. One answer is given by our main {\bf Theorem}:

\medskip \noindent
\textit{Let $\omega$ be a Delone set in Euclidean space that has
finite local complexity and uniform cluster frequencies. If $\omega$ 
is pure point diffractive, it has patch counting entropy\/ $0$.}

\medskip

A precise version is given in Theorem~\ref{main} below.  The theorem
can also be formulated for subsets of FLC Delone sets, which can be
viewed as coloured point sets, see Remark~\ref{remark-one} for
details.

In any case, the FLC
assumption is necessary in order to properly define the patch counting
function. The UCF assumption is also necessary in this theorem, as
shown by the set of visible lattice points, which satisfies FLC (as a
subset of a lattice) but not UCF (due to not being relatively dense,
and hence also not being Delone itself). In fact, the set of visible
lattice points possesses a natural autocorrelation with pure point
diffraction \cite{BMP}, while having positive patch counting entropy
\cite{Ple}.

\smallskip

Our proof of Theorem \ref{main}  relies on the following steps: 
\begin{itemize}
\item[{\rm(i)}] Pure point diffraction implies pure point dynamical
  spectrum (and vice versa).  This has recently been shown in increasing 
  degrees of generality \cite{LMS,Gou,BL,LS}.
\item[{\rm (ii)}] Pure point dynamical spectrum implies metric entropy
  $0$, as follows from a Halmos-von Neumann type representation
  theorem, see Theorem~\ref{thm:meas} below.
\item[{\rm (iii)}] There is a variational principle relating metric
  entropy and topological entropy for $\R^d$-actions, as follows from
  the work of Tagi-Zade \cite{T91}.
\item[{\rm (iv)}] The topological entropy agrees with the patch
  counting entropy for FLC Delone sets.  This is shown below in
  Theorem~\ref{HtopEqualHpatch}.
\end{itemize} 
The proof thus shows that patch counting, topological and metric
entropy coincide under the given hypotheses.  Statements (ii) and (iv)
are expected generalisations of the corresponding well-known
statements for dynamical systems with $\mathbb Z^d$-actions (compare
\cite{AL81,DGS,Pet,Wal}). They certainly belong to the mathematical
``folklore'' of the topic, see also the discussion in \cite{Robbie}.
However, no proof seems to have appeared in print so far. For this
reason, we include some details below. Let us also note that even the
various notions of entropy for an action of $\R^d$ require some care,
as (unlike in the one-dimensional case) there is no first return map.

\smallskip

Let us put our results in perspective by considering related results.
Within the regime of aperiodic order, vanishing topological entropy
has been known so far for two classes of systems that are both
characterised by their constructions. These systems are recognisable
primitive (self-affine) substitutions \cite{Rob94,HR} and uniquely
ergodic systems of finite type \cite{Wal,Shi96}.  Note that these
systems include primitive substitutions, which do \emph{not} have pure
point diffraction in general, see \cite{Q} for examples.

In this context, it is worthwhile to point out that subshifts over
finite alphabets may be strictly ergodic, but still have
positive entropy, as has been known since \cite{HK67} (see \cite{GM06}
for recent results on strictly ergodic Cantor dynamical systems with
positive entropy). In fact, there even exist almost automorphic
systems (e.g., Toeplitz flows) that are strictly ergodic and have
positive entropy \cite{BK92} (see \cite{DKL} for further results in
this direction and \cite{MP} for a generic counterpart).

On the other hand, let us mention that various random tilings are
known that quite naturally lead to FLC Delone sets that have
positive entropy (both topological and metric) as well as mixed
diffraction spectrum. They are important cases where the presence
of entropy is a clear indication of some form of disorder, which
then manifests itself also in the spectral properties of the system.

Let us briefly return to the visible lattice points mentioned above.
The natural dynamical system that arises via the closure of the
lattice translations will also contain the empty set (and various
others), and hence admits many invariant measures. However, if one
defines natural cluster frequencies via the same averaging process as
used for the natural autocorrelation, the corresponding metric entropy
vanishes \cite{Ple}. This interesting result indicates some weaker 
connection between pure point diffractiveness and vanishing metric 
entropy, which deserves further thought.

\smallskip

The paper is organised as follows. Section~\ref{Delone} is devoted to
a brief summary of Delone sets and their associated dynamical systems.
In Section \ref{PCentropy}, we discuss topological entropy and patch
counting entropy, showing their equality for FLC Delone sets in
Theorem~\ref{HtopEqualHpatch}. In Section \ref{Metric}, we recall
basic notions concerning metric entropy and its invariance under
metric isomorphism. Dynamical systems and pure point spectra are then
discussed in Section~\ref{Pure}.  In particular, we relate pure point
spectrum and entropy $0$ in Theorem~\ref{thm:meas}. Our main result,
Theorem \ref{main}, is then stated and proved in
Section~\ref{Diffraction}. In Section~\ref{Aperiodic}, we briefly
discuss how our results apply to model sets. Moreover, we provide a
result of vanishing entropy for systems whose repetitivity function
does not grow too fast. This discussion establishes vanishing entropy
for all ``standard'' examples of aperiodic order with uniform cluster
frequencies.

\section{Delone dynamical systems}\label{Delone}

We are interested in combinatorial properties of Delone sets of finite
local complexity. To such sets, dynamical systems are associated in a
natural way.  This section sets up the notation by summarising some
results from \cite{Sol,Schl,LMS}.

\smallskip
 
Consider the locally compact Abelian group $\R^d$. The 
closed ball with radius $S\geq 0$ around
$x\in \R^d$ is denoted by $B_S(x)$. We set $B_S=B_S(0)$.  The closed
hypercube with side length $n$ centred at the origin is denoted by
$C_n$.  Now, $(\varOmega,\alpha)$ is said to be a \emph{topological
  dynamical system} over $\R^d$ when $\varOmega$ is compact  and
$\alpha\! :\, \R^d\times \varOmega \longrightarrow \varOmega$ is a
continuous action of $\R^d$ on $\varOmega$. For $x\in \R^d$ and
$\omega\in\varOmega$, we write $\alpha(x,\omega)=\alpha_x(\omega)$.
The topology of $\varOmega$ gives rise to the Borel $\sigma$-algebra
$\mathcal B$ on $\varOmega$, i.e., the smallest {$\sigma$-algebra} on
$\varOmega$ which contains all open subsets of $\varOmega$.  Below, we
additionally assume that $\varOmega$ is equipped with a metric $d$
that generates the topology.

Let $m$ be an $\alpha$-invariant probability measure on
$(\varOmega,\mathcal B)$.  As usual, $m$ is called \emph{ergodic} if
the only members $B\in\mathcal B$ which are invariant up to sets of
measure $0$ satisfy either $m(B)=0$ or $m(B)=1$.  Of course, this is
equivalent to the invariant members of the completion of $\mathcal B$
all having measure $1$ or $0$ (compare Lemma~1 in
\cite[Ch.~2.1]{CFS}). The system $(\varOmega,\alpha)$ is called
\emph{uniquely ergodic} if there exists exactly one $\alpha$-invariant
probability measure on $(\varOmega,\mathcal B)$.  By standard
reasoning, see \cite{DGS,Wal}, such a measure is ergodic. When
$(\varOmega,\alpha)$ is uniquely ergodic and minimal (meaning that
each $\mathbb R^d$-orbit is dense), it is called \emph{strictly ergodic}.

\smallskip

Delone sets are discrete subsets of $\R^d$ whose points are distributed
rather uniformly.  Let positive real numbers $r$ and $R$ be given.  
A subset $\omega$ of $\R^d$ is called an $(r,R)$-\emph{Delone} set if
\begin{itemize}
\item $\omega$ is \emph{uniformly discrete} with packing radius $r$,
   i.e., any open ball of radius $r$ in $\R^d$ contains at most one point 
   of $\omega$, and
 \item $\omega$ is \emph{relatively dense} with covering radius $R$,
   i.e., any closed ball of radius $R$ in $\R^d$ contains at least one 
   point of $\omega$.
\end{itemize}
Clearly, this is only of interest when $r\le R$.
The set of all $(r,R)$-Delone sets is denoted by $\Del$.  In the
remainder of this section, we assume $0 < r \le R$ to be \emph{fixed},
with $\Del\neq \varnothing$, and usually refer to elements of $\Del$
simply as Delone sets.  The group $\R^d$ acts on $\Del$ in the obvious
way, via $\alpha_x (\omega) := -x +\omega := \, \{-x+y:y\in\omega\}$,
for each $\omega\in \Del$.  Now, let $\omega\in \Del$ be given. For
$D>0$, the elements of
\[ 
   \pos \, := \, \{(-x + \omega )\cap B_D : x \in \omega\}
\] 
are called $D$-patches of $\omega$. By construction, each $D$-patch 
contains the point 0. The cardinality of $\pos$ is
denoted by $\card(\pos)$. When $\card(\pos)<\infty$ for every $D>0$,
we say that $\omega$ has \emph{finite local complexity}, or FLC for 
short.  If, for every $D>0$ and for every $\pi \in \pos$, the limit
\[ 
    \lim_{n\to \infty}\frac{1}{|B_n|} \card \{x\in  a_n +B_n : 
     (-x + \omega) \cap B_D = \pi \}
\]
exists and is independent of the sequence $(a_n)\subset \R^d$, then
$\omega$ is said to have \emph{uniform patch frequencies}, or UCF for
short. When the corresponding property is satisfied also for all
non-empty point sets of the form $(-x+\omega)\cap K$ with
$K\subset\R^d$ compact and $x\in\omega$, which are the \emph{clusters}
of $\omega$, the point set $\omega$ has \emph{uniform cluster
  frequencies}.  Fortunately, these two notions coincide for FLC sets.
In fact, in this case, the limit can be calculated along arbitrary van
Hove sequences, see \cite{Schl} for details.

\smallskip

The notions of FLC and UCF for $\omega$ can be interpreted in terms of an
associated dynamical system. To be more precise, we need further notation. 
On $\Del$, we can introduce a metric $d$ by defining 
\[ 
  d(\xi_1,\xi_2)  := 
   \min\bigl\{ \tfrac{1}{\sqrt{2}}, \inf \{S >0\! : 
  {\exists \, u, \! v \in B_S \! :  (- u + \xi_1)\cap B_{1/S} = 
  (- v + \xi_2)\cap B_{1/S}} \} \! \bigr\}.
\]
Equipped with this metric, $\Del$ becomes a topological space on which
the action $\alpha$ of $\R^d$ is continuous, see \cite{Sol} and
references therein.  For $\omega\in \Del$, the \emph{hull} $\varOmega
(\omega)$ is defined as the closure of $\{\alpha_x (\omega) : x\in
\R^d\}$ in $\Del$.  Then, $\varOmega (\omega)$ is compact if and only
if $\omega$ is FLC.  In this case, by \cite[Thm.~3.2]{Schl},
$(\varOmega(\omega),\alpha)$ is uniquely ergodic if and only if
$\omega$ has uniform cluster frequencies, see also \cite{Sol} for an
earlier statement concerning the tiling case.

\begin{remark}
(a)
Let us mention that no metric is required for the definition of the 
hull $\varOmega (\omega)$. In fact, a more natural (and automatically
translation invariant) way works with a uniform structure, see
\cite{Schl,BL} for details. In our situation, the metric becomes
useful later on for the relation between different notions of
entropy, which is why we introduced the hull in this way. \\
(b)
The condition of relative denseness is not necessary in 
the above construction. One may furthermore consider dynamical systems
arising from suitable ensembles of uniformly discrete sets (instead 
of a single set), by taking the closure of the translation orbits
of all sets in the ensemble, such as those emerging for lattices 
gases or random tilings, compare \cite{BH} and references therein
for details. For clarity of presentation, we concentrate on the 
above setup, since we are mainly interested
in properties of FLC Delone sets.

\end{remark}

\section{Topological and patch counting entropy}\label{PCentropy}

Topological entropy for $\mathbb R^d$-actions on compact metric spaces
can essentially be defined in analogy to the case of $\mathbb
Z$-actions.  There are three equivalent approaches, via open covers,
$\ve$-separated sets, and $\ve$-dense sets, which have been studied by
Tagi-Zade \cite{T91}.  We will use a definition of entropy via
$\ve$-separated sets, as it seems most adequate for our needs.

Let $D>0$ and $\ve >0$ be given and consider $\varOmega=\varOmega(\omega)$. 
Then, a subset $\varXi$ of
$\varOmega$ is called \emph{$(D,\ve)$-separated}, if, for all
$\xi_1,\xi_2\in \varXi$ with $\xi_1\neq \xi_2$, there exists an $x\in
B_D$ with $d(\alpha_x (\xi_1),\alpha_x (\xi_2))> \ve$.

This definition can be put into more compact form.  For $D>0$ and
$\xi_1,\xi_2\in \varOmega$, define $d_D(\xi_1,\xi_2)=
\sup\{d(\alpha_x(\xi_1),\alpha_x(\xi_2)):x\in B_D\}$.  The supremum is
finite and attained since $B_D$ is compact and $\alpha$ is continuous by
assumption.  It can be shown that $d_D$ defines a metric that is equivalent
to $d$, but this is not needed in the sequel. Now, $\varXi\subset\varOmega$
is $(D,\ve)$-separated if and only if, for all $\xi_1,\xi_2\in \varXi$ with
$\xi_1\neq \xi_2$, we have $d_D(\xi_1,\xi_2)>\ve$.

\smallskip

By compactness of $\varOmega$, there exists an upper bound on the
number of elements of any $(D,\ve)$-separated set. So,
$ N(D,\ve)  :=  \max\{\mbox{card}(\varXi) : \varXi \subset \varOmega
     \,\mbox{ is $(D,\ve)$-separated}\}$ exists and is finite. Define now
$H_{\ve}  :=  H_{\ve} (\varOmega,\alpha) := \limsup_{n\to \infty} \frac{1}{|B_n|}
    \log N(n,\ve) \, \ge \, 0 $ as usual.
Obviously, $H_{\ve}\geq H_{\ve'}$ whenever $\ve \leq \ve'$. 
Thus, the limit
\[ 
   h_{\rm top}  \, := \, h_{\rm top} (\varOmega,\alpha) \, := \,
    \lim_{\ve \to 0} H_{\ve} (\varOmega,\alpha)
\]
exists in $\R_{\ge 0}\cup\{\infty\}$. 
It is called the \emph{topological entropy} of $(\varOmega,\alpha)$.  
By compactness of $\varOmega$, equivalent metrics on $\varOmega$ lead to
the same topological entropy, see also \cite[Thm.~7.4]{Wal}. In fact, it
only depends on the topology, and can be defined independently of the
metric, compare \cite{T91,DGS}. It is thus clear that topological entropy
is invariant under topological conjugacy. In view of our later
calculations, we prefer to use a formulation with a given metric.

Of course, various normalisations are possible. For example, one can
replace the balls $B_n$ by hypercubes $C_n$.  This may lead to
some overall factor in the definition of the entropy, as a consequence
of the volume ratio of balls versus cubes of the same extension.
However, this cannot
influence whether the entropy is $0$ or not. In the cases we
deal with, the entropy will turn out to be $0$, and such factors will
thus not matter at all.

\smallskip

Recall that the cardinality of the set $\pos$ of all $D$-patches is
denoted by $\mbox{card}(\pos)$. The \emph{patch counting entropy}
$h_{\rm pc}(\omega)$ of $\omega$ is defined by
\begin{equation} 
   h_{\rm pc}(\omega) \, := \,
   \limsup_{n\to \infty} \frac{1}{|B_n|} 
   \log \bigl(\mbox{card}(p_\omega(n))\bigr) \/ ,
\end{equation} 
which exists and is finite, see \cite[Thm.~2.3]{Lag}.  The following
theorem states that the patch counting entropy of an FLC Delone set
coincides with the topological entropy of the associated dynamical
system. This is an extension of a well-known result in symbolic
dynamics, compare \cite[Ch.~6.3]{Pet}.  In order to prepare its proof,
we first derive a lower bound for the distance of two Delone sets that
share the origin as a vertex.

\begin{lemma}\label{geometry}
Let $\{\xi_1,\xi_2\}\subset\Del$ satisfy $0\in\xi_1\cap\xi_2$ and 
$\xi_1\neq \xi_2$. Then, if $S>0$ is any number that satisfies 
$\xi_1\cap B_S\ne \xi_2\cap B_S$, one has
$d(\xi_1,\xi_2) \ge  \min\bigl\{\frac{1}{\sqrt{2}}, \frac{r}{2},
   \frac{1}{S}\bigr\}$.
\end{lemma}

\begin{proof} 
  Note that $\xi_1\ne\xi_2$ by assumption.  Fix $S>0$ such that
  $\xi_1\cap B_S\ne \xi_2\cap B_S$. Choose $\widetilde r\le
  \min\{r/2,1/S\}$.  Then, for $x,y\in B_{\widetilde r}$ arbitrary, we
  have $(\xi_1-x)\cap B_{1/\widetilde r}\ne (\xi_2-y)\cap
  B_{1/\widetilde r}$, as can easily be seen by considering the cases
  $x=y$ and $x\neq y$ separately. This implies $d(\xi_1,\xi_2)\ge
  \min\{\widetilde r, 1/\sqrt{2}\}$, and the claim follows.
\end{proof}

Note that the above proof does not require relative denseness
of $\xi_1$ or $\xi_2$.

\begin{theorem} \label{HtopEqualHpatch} 
Let $\omega\in \Del$ be an FLC Delone set, with hull $\varOmega(\omega)$. 
Then, for each $\ve<\ve_0  :=  
\min\bigl\{\frac{1}{\sqrt{2}},\frac{r}{2},
\frac{1}{2R}\bigr\}$,
one has $H_{\ve} (\varOmega(\omega),\alpha) = 
h_{\rm pc}(\omega)$ . In particular,
\[
   h_{\rm pc}(\omega) \, = \, h_{\rm top} (\varOmega(\omega),\alpha).
\]
\end{theorem}

\begin{proof} {}For $D>0$, we use the abbreviation
  $n(D):=\mbox{card}(p_\omega(D))$.  Now, consider a subset
  $\varXi^{(D)}=\{\xi_1,\ldots,\xi_{n(D)}\}$ of $\varOmega(\omega)$
  that represents the $D$-patches around the point $0\in \R^d$, i.e.,
  $\varXi^{(D)}$ has $n(D)$ elements and $\pos = \{ \xi \cap B_D : \xi
  \in \varXi^{(D)}\}$.  We show two inequalities.

\smallskip

``$\leq$'': 
Assume $n(D)>1$ and fix some $1\le i,j\le n(D)$ 
with $i\ne j$. Choose a number $0<S_{ij}\le D$ such that 
$\xi_i\cap B_{S_{ij}}\ne\xi_j\cap B_{S_{ij}}$ and $\xi_i
\cap B_{S}=\xi_j\cap B_{S}$ for all $S<S_{ij}$. Distinguish 
two cases. If $S_{ij}\le 2R$, Lemma~\ref{geometry} implies
\[
  d_D(\xi_i,\xi_j) \, \ge \, d(\xi_i,\xi_j)
  \, \ge \, \min\bigl\{  \tfrac{1}{\sqrt{2}}, \tfrac{r}{2},
   \tfrac{1}{S_{ij}}\bigr\} \, \ge \, \min\bigl\{\tfrac{1}{\sqrt{2}},
   \tfrac{r}{2},\tfrac{1}{2R}\bigr\} \, = \, \ve_0 \/ .
\]
If $S_{ij}>2R$, fix $x_0\in \partial B_{S_{ij}}$ such that 
$x_0\in \xi_i\cup\xi_j$ but $x_0\notin \xi_i\cap\xi_j$. Choose 
a ball $B\subset S_{ij}$ of radius $R$ such that $x_0\in B$. 
Choose $x\in B\cap \xi_i$. Such a point exists, since $\xi_i$ 
is relatively dense with radius $R$. By construction, we also 
have $x\in \xi_j$. Since $x_0\in B_{2R}(x)$, we have
 $\xi_i\cap B_{2R}(x) \ne  \xi_j\cap B_{2R}(x)$.
We can now apply Lemma~\ref{geometry}, with $x\in B_D$ considered
to be the origin. This gives
\[
  d_D(\xi_i,\xi_j)\, \ge\, d(\alpha_{-x}(\xi_i),\alpha_{-x}(\xi_j))
  \,\ge\, 
   \min\bigl\{\tfrac{1}{\sqrt{2}}, \tfrac{r}{2},\tfrac{1}{2R}\bigr\}
   \,\ge\, \ve_0 \/ .
\]
Combining the two cases, we infer that, for any $D>0$ and any
$\ve<\ve_0$, the set  $\varXi^{(D)}$ is indeed 
$(\ve,D)$-separated.
(If $n(D)=1$, this is trivially the case.) This implies
\[ 
  n(D) \, = \, \mbox{card}(\varXi^{(D)}) \,\leq\,  N(D,\ve)\/ ,
\] 
and $h_{\rm pc}(\omega) \leq H_{\ve}$ follows. As $\ve$ was arbitrary 
with $0  <\ve < \ve_0$, we obtain ``$\leq$''. 

\smallskip

``$\geq$'': 
Let $\ve >0$ and $D>0$ be arbitrary. Define $\rho(D):=
D + R +\ve+1/\ve$. Let $\varXi$  be an arbitrary
$(D,\ve)$-separated set. For any $\xi\in \varXi$, we can find an
$x_\xi\in\xi\cap B_R$, since $\xi$ is relatively dense with radius $R$.
Thus, we can define a mapping
$ \varPsi \!:\, \varXi \longrightarrow \varXi^{(\rho(D))}$,
such that $\xi\in \varXi$ is mapped to $\eta \in  \varXi^{(\rho(D))}$ 
whenever  $ (-x_\xi + \xi)\cap B_{\rho(D)}  = \eta \cap B_{\rho(D)}$.
In general, the map $\varPsi$ is not injective. A bound
on the number of preimages of $\eta$, meaning $\card(\varPsi^{-1}(\eta))$, 
is obtained as follows. First, we claim that 
\begin{equation}\label{distance}
   d(x_{\xi_1},x_{\xi_2}) \, \geq \, \ve
\end{equation}
whenever $\xi_1 \neq \xi_2$ are elements of $\varXi$ with $\varPsi(\xi_1)
= \varPsi (\xi_2)$. To see this, assume the contrary. Then,
a short calculation yields $d_D(\xi_1,\xi_2)< \ve$, based on the definition 
of $\rho(D)$ and the construction of $\varPsi$. However,
this is impossible, as $\xi_1$ and $\xi_2$ are different elements
of a $(D,\ve)$-separated set $\varXi$.  Furthermore, by
compactness of $B_R$, there exist $M=M(\ve)\in \N$ and
$x_1,\ldots, x_M\in B_R$ such that 
\begin{equation}\label{covering}
   B_R \,\subset\, \bigcup_{i=1}^M B_{\ve/2}(x_i).
\end{equation}
Combining \eqref{distance} with \eqref{covering}, we 
infer that no $\eta\in\varXi^{(\rho(D))}$ can have 
more than $M$ inverse images under $\varPsi$. This implies
$\mbox{card}(\varXi) \le M(\ve) \cdot n(\rho(D))$.
As $\varXi$ is an arbitrary  $(D,\ve)$-separated set, 
we infer  $N(D,\ve) \le M(\ve)\cdot n(\rho(D))$.
With $\lim_{D\to \infty} \frac{1}{D} \rho(D) = 1$, one
finds $H_{\ve} \leq h_{\rm pc}(\omega)$. As $\ve>0$ was arbitrary, 
the desired inequality follows, and the claimed equality is shown.
\end{proof}

\begin{remark} \label{remark-one} 
  {\rm (a)} The above proof can be adapted to coloured FLC Delone sets
  (with finitely many colours), compare \cite{LMS,LStoll} for a
  discussion of coloured Delone sets.  In this case, the FLC property
  is inherited, but the UCF property becomes increasingly more
  restrictive. Every subset of an FLC Delone set gives rise to a
  coloured FLC Delone set, via the obvious colouring with two colours.
  An example is given by the visible lattice points. Note that the
  number of $D$-patches (and possibly also the corresponding
  entropies) of an FLC Delone subset and its coloured counterpart may
  be different.
  \\
  {\rm (b)} The above proof uses the Delone property of the underlying
  set in both directions. It is an open question whether the condition
  of relative denseness is necessary for the validity of the theorem.
  See Example~\ref{bob} below for an FLC set which is
  not Delone. \\
  {\rm (c)} An alternative proof can be given, using the definition of
  topological entropy by open covers, in analogy to the case of
  $\mathbb Z$-actions.  Note, however, that the dynamical system is
  not expansive, and entropy is computed by a sequence of open covers
  with diameters
  that shrink to $0$.\\
  {\rm (d)} The above result generalises the cases of $\mathbb
  Z$-actions \cite[Ch.~6.2]{Pet} and $\mathbb Z^d$-actions \cite{AL81}.
\end{remark}

\begin{ex}
  Quasiperiodic tilings provide, via a natural point decoration,
  interesting examples of aperiodic FLC Delone sets. This includes the
  majority of tilings used in crystallography and physics to study the
  properties of quasicrystals. Examples are, among various others, the
  Penrose tiling, the Ammann-Beenker tiling, the shield tiling, the
  Ammann-Kramer tiling, and the Danzer tiling. See the review \cite{B}
  and references therein. The underlying point sets are regular model
  sets, and some tilings are substitution tilings. See Section
  \ref{Aperiodic} for a discussion of their entropy and diffraction
  properties.
\end{ex}

\begin{ex} \label{rantil} 
  It is well-known that fully-packed dimer models on lattices or
  periodic graphs give rise to interesting random tiling examples
  \cite{Ric1,Ric2} and, via a natural point decoration, also to FLC
  Delone sets. A random tiling ensemble gives rise to a dynamical
  system, via the closure of the translations of the Delone sets
  derived from all tilings of the ensemble. The system has positive
  patch counting (and hence topological) entropy, and a diffraction
  spectrum of mixed type, see \cite{BH} for a proof of this statement.
  For the dimer model on the square lattice (the domino case) and on
  the hexagonal lattice (the lozenge case), the corresponding
  equilibrium (or Gibbs) measure is unique, so that also the metric
  entropy is positive.  The topological entropy of both models is of
  the form
\[
   m(P)\, = \, \int_0^1\int_0^1\log|P(e(s),e(t))| \,{\rm d}s\,{\rm d}t\, ,
\]
for certain two-variable Laurent polynomials $P=P(x,y)$, where
$e(t)=\exp(2\pi {\rm i} t)$. We have $P(x,y)=4+x+1/x+y+1/y$ in the
domino case \cite[Eq.~8]{K}, and $P(x,y)=1+x+y$ in the lozenge case
\cite[Eq.~5]{Wu}.

The quantity $m(P)$, the average of $\log|P|$ over the real 2-torus,
is also called the (logarithmic) \emph{Mahler measure} of $P$. For $P$
having integer coefficients only, it arises as an entropy within a
certain $\mathbb Z^2$-dynamical system associated to $P$, see
\cite{LSW}.  It would be interesting to investigate its relation to
the $\mathbb R^2$-dynamical system discussed above, see also
\cite{So}.

\end{ex}

\begin{ex} \label{bob}
  Let us briefly discuss a situation that sheds some light on the
  connection between FLC sets and Delone sets\footnote{We thank Robert
    V.\ Moody for suggesting this type of example.}. Consider Euler's
  number $e$ and define $a^{}_{n} = 1 + e + e^2 + \ldots + e^{n-1}$
  for $n\in\N$.  The set $\varLambda=\{\pm a^{}_{n} : n\in\N\}$ is
  uniformly discrete, but not relatively dense, and hence not Delone.
  Since $a^{}_{n+1} - a^{}_{n}=e^{n} \xrightarrow{n\to\infty}\infty$,
  there are only \emph{finitely} many patches of a given diameter, so
  that $\varLambda$ is an FLC set. As the patch counting function
  grows only logarithmically with patch size, one can conclude that
  $h_{\rm pc}=0$ for $\varLambda$.

  On the other hand, as $e$ is transcendental, the points of
  $\varLambda$ are rationally independent, so that
  $[\varLambda-\varLambda]$ is not a finitely generated $\Z$-module.
  Consequently, by a result from \cite{Lag}, $\varLambda$ cannot be a
  subset of any FLC Delone set. Since $\varLambda$ has vanishing density,
  its entropy per volume does not seem to be the right quantity
  to look at. However, when reconsidering it as a two-colour
  Delone set (which is then \emph{not} FLC), it is not clear how
  to relate the various entropies, or what their values are.
\end{ex}

\section{Metric entropy and variational principle} \label{Metric} 

In this section, we discuss the metric entropy, its invariance under
metric isomorphism and the variational principle.  While we need two
properties of the metric entropy, the details of its actual definition
play no role below. Therefore, we omit the technicalities and
concentrate on these properties instead. For details concerning its
definition and the variational principle, we refer the reader to the
paper of Tagi-Zade \cite{T91}.

\smallskip

The first part of our discussion does not require a topological
dynamical system.  It works whenever we are given a measure dynamical
system, i.e., a measurable action $\alpha$ of $\R^d$ on the measurable
space $\varOmega$ together with an $\alpha$-invariant probability
measure $m$. To these data, one can associate a quantity called 
\emph{metric entropy} or Kolmogorov-Sinai entropy, see \cite{CFS,DGS}
for details.

Let us continue by explaining the notion of metric isomorphism of
measure preserving transformations. Let two measure spaces
$(\varOmega_i, \mathcal B_i, m_i)$, $i\in\{1,2\}$, be given. Let
$\alpha_i$ be measure preserving actions of $\R^d$ on $\varOmega_i$.
Recall that the measure algebras $(\widetilde{\mathcal B}_i,\widetilde
m_i)$ arise from the corresponding $\sigma$-algebras by identifying
elements whose symmetric difference has measure $0$.  We then say that
$(\varOmega^{}_1, \alpha^{}_1)$ and $(\varOmega^{}_2,\alpha^{}_2)$ are
\emph{metrically isomorphic} (or conjugate modulo $0$) if there is a
measure algebra isomorphism
\[
   \varPhi \! : \; (\widetilde{\mathcal B}_2,\widetilde m_2) 
   \longrightarrow (\widetilde{\mathcal B}_1,\widetilde m_1)\, ,
   \quad \mbox{with}\;\: \varPhi(
   \alpha_{1,x} (V) ) = \alpha_{2,x} (\varPhi (V)) \/ ,
\]
for all $x\in \R^d$. A direct consequence of the definitions of metric
entropy and metric isomorphism is the 
following well-known fact, compare \cite{CFS,DGS,Wal}.

\begin{lemma}\label{lem:invariant}
   Metric entropy is invariant under metric isomorphisms.  \qed
\end{lemma}

Let us now turn to the variational principle, which relates
metric and topological entropy.  

\begin{theorem} \label{var-princ}
  Let $\varOmega$ be a compact metric space, equipped with the Borel
  $\sigma$-algebra and a continuous action $\alpha$ of\/ $\R^d$. Then,
  the topological entropy is the supremum of the metric entropies 
  taken over all $\alpha$-invariant probability measures on $\varOmega$.
\end{theorem}

This is proved in \cite{T91} by first relating the 
$\mathbb R^d$-entropies to corresponding $\mathbb Z^d$-entropies 
via restricting the underlying translation group. Since these 
entropies are shown to coincide, the statement follows from the 
variational principle for $\mathbb Z^d$-actions, which is
discussed in \cite[Sec.~6]{Rue}. It is an extension of the
case $d=1$ from \cite{DGS,Wal}, see also \cite{El} and references
given there for further details.

\section{Pure point spectrum and the representation theorem}
\label{Pure} 

{}First, we recall basic notions concerning pure point spectra and
formulate the representation theorem. This circle of ideas is
sometimes discussed under the name Halmos-von Neumann theorem, see
\cite{Wal,CFS,Sinai} for details, and \cite{Robbie} for its recent
appearance in a related context. We then present an abstract result
relating pure point spectrum and vanishing metric entropy for ergodic
dynamical systems on compact metric spaces.

\smallskip

A measurable dynamical system $(\varOmega,\alpha,m)$ gives rise to a
unitary representation $T = T^{(m)}$ of $\R^d$ on $\Lo$, which is
given by
\[ 
   T_x f = f \circ \alpha_x \quad
        \mbox{for  $f\in \Lo$ and $x\in \R^d$}.
\]
An $f\in \Lo$ is called an \emph{eigenfunction} of $T$ if $f\neq 0$
and if there exists a $\lambda$ in the dual group $\widehat{\R^d}$ of
$\R^d$ with $T_x f= (\lambda,x) f$ for all $x\in \R^d$. If $\Lo$ has a
basis that purely consists of eigenfunctions of $T$, $T$ is said to
have \emph{pure point dynamical spectrum}. Alternatively, one then
also says that $(\varOmega,\alpha,m)$ has pure point spectrum.

There is a representation theorem for ergodic measure preserving
$\mathbb R^d$-actions with pure point spectrum, as for ergodic
$\mathbb Z$-actions. It can be proved by adapting the arguments for
$\mathbb Z$-actions, as given in \cite[Thms.~3.4 and 3.6]{Wal}, to the
case of $\mathbb R^d$-actions, see also the detailed exposition in
\cite[Ch.~I.5]{Sinai} for the case $d=1$. For the convenience of the
reader, we provide some details for the general case.

An action $\beta$ of $\R^d$ on a compact group $\T$ 
is called a \emph{rotation} if there exists a group homomorphism
$i \! : \; \R^d \longrightarrow \T$ such that 
$\beta_x (\xi) =  i(x) \xi$, for all $x\in \R^d$.
If $i$ has dense range, then $(\T,\beta)$ is uniquely ergodic, 
as the normalised Haar measure is the unique translation invariant 
probability measure on $\T$. From now on, 
$\alpha$ denotes an $\R^d$-action.

\begin{theorem}\label{thm:rep}
  Let $(\varOmega,\alpha, m)$ be an ergodic measurable dynamical
  system with pure point dynamical spectrum. Then, $\alpha$ is
  metrically isomorphic to a uniquely ergodic rotation $\beta$ 
  on some compact Abelian group $\mathbb T$ (equipped with its 
  unique normalised Haar measure). If
  $L^2(\varOmega,m)$ admits a countable orthonormal basis, $\T$
  is metrisable with a $\beta$-invariant metric.
\end{theorem}

\begin{proof}
Denote the $\alpha$-invariant probability measure on 
$\varOmega$ by $m$. By assumption, there exists $P\subset
\widehat{\R^d}$ and a basis $ \{ f_\lambda\}_{\lambda\in P}$ on 
$\Lo$ such that $f_\lambda$ is an eigenfunction of $T$ to 
$\lambda\in P$. (Here, we use that every eigenvalue has 
multiplicity 1 due to ergodicity.) 

Now, equip $P$ with the discrete topology and denote its dual group by
$\T$. Then, $\T$ is a compact Abelian group whose dual is $P$.  As
$P\subset \widehat{\R^d}$, the natural map $i\! :\,
\R^d\longrightarrow \T$, defined by $i(t)(\lambda) :=(\lambda,t)$, has
dense range.  It induces a continuous action
\[ 
   \beta \! :\, \R^d\times \T \longrightarrow \T \, ,
   \quad (x,\sigma) \mapsto i(x)\sigma \/ .
\]
As $i$ has dense range and the Haar measure is unique, this 
action is uniquely ergodic. We now show that $(\varOmega,\alpha, m)$ 
is metrically isomorphic to $(\T,\beta)$. This proceeds in three steps, 
see also \cite[Thm.~3.4]{Wal}: 

Step 1: Without loss of generality, we may assume that $f_\lambda 
f_\mu = f_{\lambda + \mu}$, $\forall \lambda, \mu\in P$. 

\textit{Argument:} Follow the proof of \cite[Thm.~3.4]{Wal} or
   \cite[Ch.~I.5]{Sinai}.

\smallskip

Step 2: There exists a map $M \!:\, L^2 (\varOmega,m)\longrightarrow 
L^2(\T )$ which commutes with the action of $\R^d$ and satisfies 
$M(f g) = M(f) M(g)$ and $M (\bar{g}) = \overline{M(g)}$, 
$\forall f,g\in L^\infty (\varOmega,m)$. 

\textit{Argument:} Define $j\! :\, \R^d\longrightarrow \widehat{\T}$ 
by $j(x) (\rho) :=(\rho,  i(x))$. Thus, each $j(x)$ is a 
character and in fact a bounded eigenfunction on $\T$. Then,
the map $\sum c_\lambda f_\lambda \mapsto \sum c_\lambda\,
j(\lambda)$ has the desired properties by Step 1. 

\smallskip

Step 3: There exists an isomorphism of measure algebras, 
which commutes with the action of $\R^d$. 

\textit{Argument:} Apply $M$ to characteristic functions. 

\smallskip

Finally, we discuss metrisability: If $L^2(\varOmega,m)$ permits a
countable orthonormal basis, then so does $\T$. Thus, the set of
characters of $\T$ (which gives an orthonormal basis) is countable.
Let $\chi_n, n\in \N$, be an enumeration of the characters and define
a metric by
\[
   d_{\mathbb T}(\xi,\eta) \, := \, \sum_{n\in \N} 2^{-n} 
   \lvert (\chi_n, \xi - \eta)\rvert \, . 
\]
Then, $d_{\mathbb T}$ is invariant as each character is an
eigenfunction of $\beta$.  Obviously $d_{\mathbb T}$ is also
continuous, and it separates points by Pontryagin duality. Thus,
it generates the topology of the compact group $\T$.
\end{proof}

\begin{remark}
  When the Delone set $\omega$ is a model set, see \cite{Moody} for
  details, it emerges from a cut and project scheme, which provides a
  natural candidate for the compact group $\T$ as the factor group of
  the embedding space by the embedding lattice. This can be made
  explicit by means of the \emph{torus parametrisation}, see
  \cite{BHP,HRB,Schl,BLM} for details.
\end{remark}

\begin{lemma} \label{lem:$0$} 
  Let $(\varOmega,\alpha)$ be a topological dynamical system with an
  $\R^d$-invariant metric $d$ on $\varOmega$. Then, one has
  $h_{\rm top} (\varOmega,\alpha)=0$.
\end{lemma}
\begin{proof} 
As the metric $d$ is $\R^d$-invariant, $N(D,\ve)$ does not depend on $D$.  
Thus, $H_{\ve}=0$ for every $\ve>0$, and the statement follows.
\end{proof}

\begin{theorem}\label{thm:meas}
  Let $(\varOmega,\alpha, m)$ be a measurable dynamical system
  that is ergodic and has a pure point dynamical spectrum, with
  a countable basis of eigenfunctions.  
  Then, $(\varOmega,\alpha,m)$ has metric entropy $0$.
\end{theorem}

\begin{proof} 
  By Theorem \ref{thm:rep}, $(\varOmega,\alpha,m)$ is metrically
  isomorphic to a rotation $(\mathbb T,\beta)$, which is 
  uniquely ergodic by Theorem~\ref{thm:rep}, with Haar measure
  $\mu$ as the invariant measure.
  As we have a countable basis, $\T$ admits a $\beta$-invariant metric.  
  Lemma~\ref{lem:$0$} then gives
  that $(\mathbb T,\beta)$ has topological entropy $0$.  Since
  topological entropy is an upper bound for measure theoretic entropy
  (where we actually get equality here due to unique ergodicity),
  $(\mathbb T,\beta,\mu)$ has metric entropy $0$. 
  Since $(\mathbb T,\beta,\mu)$ and $(\varOmega,\alpha,m)$
  are metrically isomorphic, Lemma~\ref{lem:invariant} implies that
  $(\varOmega,\alpha,m)$ has metric entropy $0$ as well.
\end{proof}

\section{Pure point diffraction and entropy $0$} \label{Diffraction}

In this section, we finally study the connection between pure point
diffraction and entropy. We start by discussing diffraction, as
initiated in a mathematical setting by Hof \cite{Hof}, compare
\cite{Schl,Gou,BL,LS} and references given there for recent results.
Let $\omega\in \Del$ with hull $\varOmega = \varOmega(\omega)$ be
given.  Let $m$ be a translation invariant measure on the hull. Then,
$m$ induces a measure, called \emph{autocorrelation} of $\varOmega$
and denoted by
$\gamma^{}_{\varOmega,m} =\gamma^{}_\varOmega$, on $\R^d$ as follows:
Choose a continuous function $\sigma$ with compact support and
$\int_{\R^d} \sigma (t) {\rm d}t =1$. For $\varphi$ continuous with
compact support on $\R^d$, we then set
\[
  \gamma^{}_\varOmega (\varphi ) \, = \,
  \int_\varOmega \Bigl(\,\sum_{s,t\in \omega} 
  \sigma (s) \varphi (s-t) \Bigr) \,\mathrm{d} m(\omega).
\]
It turns out that the measure $\gamma^{}_{\varOmega}$ does not depend
on the choice of $\sigma$. It is not hard to see that
$\gamma^{}_\varOmega$ is positive definite.  Thus, its Fourier
transform $\widehat{\gamma^{}_\varOmega}$ exists and is a positive
measure. It is called the \emph{diffraction measure} of $\varOmega$.
When the Fourier transform $\widehat{\gamma^{}_\varOmega}$ is a pure
point measure, we say that $\varOmega$ is \emph{pure point
  diffractive}.

In the case of finite local complexity and uniform cluster frequencies,
the autocorrelation measure $\gamma^{}_{\omega}$ of any fixed element 
$\omega\in\varOmega$ can be computed as
\[
   \gamma^{}_\omega \, = \, 
   \lim_{n\to \infty} \frac{1}{\lvert B_n \rvert} 
   \sum_{x,y\in \omega\cap B_n} \delta_{x-y} \, ,
\]
where $\delta_z$ denotes the unit point mass 
at $z\in \R^d$, and the limit is taken in the vague topology. 
In fact, the limit can be computed along arbitrary 
van Hove sequences, and does not depend on the choice of
$\omega$, due to unique ergodicity. So, $\gamma^{}_{\omega}
=\gamma^{}_{\varOmega}$ for all $\omega\in\varOmega$ in
this case. In general, whenever $\gamma^{}_{\omega}$ is well-defined
(possibly with respect to a specified averaging sequence),
its Fourier transform exists and is a positive measure. One calls 
the single Delone set $\omega$ \emph{pure point diffractive} when 
$\widehat{\gamma^{}_{\omega}}$ is a pure point measure.

\smallskip
The following result from \cite{LMS} (see \cite{Gou,BL,LS} for
further generalisations) gives a crucial characterisation of pure
point diffractiveness. 

\begin{lemma} \label{characterization} 
  Let $\omega\in \Del $ be given, with finite local complexity, 
  hull $\varOmega=\varOmega(\omega)$ and invariant probability
  measure $m$.  Then, the following assertions are equivalent:
\begin{itemize}
\item[(i)] $\varOmega$ is pure point diffractive, i.e., 
   $\widehat{\gamma^{}_\varOmega}$ is a pure point measure. \smallskip
\item[(ii)] $(\varOmega,\alpha,m)$ has pure point dynamical spectrum.
   \qed
\end{itemize}
\end{lemma}

This lemma and the considerations of the preceding sections imply the 
following.

\begin{coro} \label{coro-one}
  Let $\omega\in \Del $ be given, with finite local
  complexity, hull $\varOmega=\varOmega(\omega)$ and ergodic measure $m$. 
  If $\widehat{\gamma^{}_\varOmega}$ is pure point, 
  $(\varOmega,\alpha,m)$ has metric entropy\/ $0$.
\end{coro}
\begin{proof} 
  By Lemma~\ref{characterization}, $(\varOmega,\alpha,m)$ has pure point
  dynamical spectrum. Consequently, the metric entropy vanishes by
  Theorem~\ref{thm:meas}.
\end{proof}

We can now state and prove our main result. 

\begin{theorem}\label{main}  
  Let $\omega\in \Del$ be given, with finite local
  complexity and uniform cluster frequencies.  If $\omega$ is
  pure point diffractive, one has $h_{\rm pc}(\omega)=0$.
\end{theorem}

\begin{proof}  
  As $\omega$ has uniform cluster frequencies,
  $(\varOmega(\omega),\alpha)$ is uniquely ergodic by
  \cite[Thm.~3.2]{Schl}. Denote the unique $\alpha$-invariant ergodic
  measure by $m$. Corollary~\ref{coro-one} and the variational
  principle from Theorem~\ref{var-princ} then imply that $h_{\rm top}
  (\varOmega(\omega),\alpha) =0$.  By Theorem \ref{HtopEqualHpatch},
  we obtain that the patch counting entropy satisfies
\[
    h_{\rm pc}(\omega) \, = \, h_{\rm top} (\varOmega(\omega),\alpha) 
      \, = \, 0 \, ,
\]
which completes the argument.
\end{proof}

\begin{remark} 
  As we have not assumed repetitivity of $\omega$, the hull
  $\varOmega(\omega)$ may be bigger than the LI-class of $\omega$ (so
  we need not have strict ergodicity).  Consequently, different
  elements of $\varOmega(\omega)$ may possess different patch counting
  functions. However, by the structure of the hull, they are all
  majorised by the patch counting function of $\omega$ itself (note
  that one can only ``lose'' patches in the limiting process).
  Consequently, with $h_{\rm pc}(\omega) = 0$ from Theorem~\ref{main},
  one also has $h_{\rm pc} (\omega')=0$ for all
  $\omega'\in\varOmega(\omega)$.
\end{remark}

\section{Aperiodic order and entropy $0$}\label{Aperiodic}

There are two main classes of examples of aperiodic order. One class
consists of model sets, and the other is given by primitive
substitutions. In this section, we discuss the vanishing of the
entropy for these classes.

\smallskip

As regular model sets have uniform cluster frequencies and are pure 
point diffractive \cite{Schl}, our main result immediately implies:

\begin{theorem} 
   Regular model sets have patch counting and topological 
   entropy $0$.   \qed
\end{theorem}

{}For primitive self-affine substitutions, patch counting entropy $0$
is shown in \cite{HR}. Here, we discuss a simple criterion for entropy
$0$ whenever the repetitivity function does not grow too fast. In
particular, this criterion applies to linearly repetitive systems as
introduced in \cite{Du,LP1}, which include primitive self-similar
substitutions \cite{LP1, Sol2}. We thus recover a part of the result
from \cite{HR} by a different method.

\begin{definition} 
  A Delone set $\omega$ has \emph{repetitivity function} 
  $F \! : \,[1,\infty) \longrightarrow [1,\infty)$ if any 
  patch of size $D\ge 1$ is contained in any patch of size 
  $F(D)$, i.e., if for any $x\in \omega$ the inclusion
\[
   p_\omega (D) \, \subset \,
   \{ (\omega - y)\cap B_D : y \in \omega \cap B_{F(D)} (x)\}
\]
holds. If $F$ can be chosen to be $F(D) = c D$ for some $c>1$, $\omega$ 
is called \emph{linearly repetitive}. 
\end{definition}

\begin{theorem} \label{thm-seven} 
  If $\omega$ is a Delone set in $\R^d$ with repetitivity function
  $F$, there exists some $\kappa >0$ with $\mbox{\rm card}
  (p_\omega (D)) \leq \kappa\, \lvert B_{F(D)} \rvert $ for all $D\geq
  1$. Thus, $h_{\rm pc}(\omega)=0$ whenever one has $F(D) = o(\exp (D))$.  
  In particular, $h_{\rm pc}(\omega) = 0$ when $\omega$ is linearly
  repetitive.
\end{theorem}

\begin{proof} 
By the Delone property, there exists a constant $\kappa_1>0$ such that
\[
  \mbox{card}(\omega \cap B_D (p)) \,\leq\, \kappa_1 \lvert B_D \rvert
\]
for all $D\geq 1$ and $p\in\R^d$. Let $D\geq 1$ be arbitrary and choose 
some $x\in \omega$. By definition of the repetitivity function, 
every patch of size $D$ can be found in $B_{ F(D)} (x)$. Obviously, 
there cannot be more patches of size $D$ in $B_{F(D)} (x)$ than there 
are points of $\omega$ in $B_{F(D)}$. This gives
  $  \mbox{card}(p_\omega (D)) \leq 
    \mbox{card}(\omega\cap B_{F(D)} (x))
    \, \leq \, \kappa_1 \vert B_{F(D)}\rvert$, 
which proves the statement. 
\end{proof}

\begin{remark} (a) Note that the proof of Theorem~\ref{thm-seven} 
does not use relative denseness. \\
(b) For aperiodic linearly repetitive sets, there is also a
corresponding lower bound. This had been conjectured by Lagarias and
Pleasants \cite{LP1} and was later proved in \cite{Len}. In
particular, for primitive self-similar substitutions, the patch
counting function grows asymptotically proportional to $|B_D|$, see
also the discussion in \cite{Ric}. \\
(c) Linear repetitivity (and thus $h_{\rm pc}(\omega) = 0$) alone
does \emph{not} imply the absence of continuous diffraction components,
as can be seen from the Thue-Morse chain (singular continuous parts)
or the Rudin-Shapiro chain (absolutely continuous parts),
see \cite{Q,HB} for details.

\end{remark}

\bigskip

\section*{Acknowledgements}

It is our pleasure to thank J.~Bellissard, R.V.~Moody and A.C.D.\ van
Enter for valuable discussions and various hints on the history of the
subject. We are particularly grateful to J.~Kwapisz for helpful
suggestions concerning Theorem~\ref{thm:meas} and the variational
principle, and to B.~Weiss for pointing out reference \cite{T91} to
us.  We thank the referees for several useful hints that helped to
improve the presentation and the connection with other results.  This
work was supported by the German Research Council (DFG), via an
individual project (DL) and within the CRC 701.

\end{document}